\documentclass[10pt,twoside]{article}
\usepackage{amsfonts,amsmath,amssymb,amsopn,amsthm,mdwlist}
\usepackage[cns,uns]{llconnex}

\input diagxy

\newcommand{\op}{{\rm op}}

\newcommand{\ladj}{\dashv}

\newcommand{\too}{\longrightarrow}

\newcommand{\oto}[1]{\stackrel{\underrightarrow{\stackrel{#1}{\,\,\,\,\hphantom{\too}}}}{}}  

\theoremstyle{definition}

\author{\sc Brian Day\footnote{Department of Mathematics, Macquarie University, NSW 2109, Australia.}}
\date{February 13, 2009}
\title{When is Existential Quantification Conservative?}

\begin{document}

\maketitle

\begin{abstract}
We describe a sufficient condition for the process of left Kan extension to be a conservative functor. This is useful in the study
of graphic Fourier transforms and quantum categories and groupoids.
\end{abstract}

\section{Introduction}

\newcommand{\V}{\mathcal V}
\newcommand{\power}{\mathcal P}
\newcommand{\A}{\mathcal A}
\newcommand{\C}{\mathcal C}
\newcommand{\E}{\mathcal E}
\newcommand{\Lan}{\rm Lan}
\newcommand{\lbrak}{\left[}
\newcommand{\rbrak}{\right[}
Let $\V$ be a complete and cocomplete symmetric monoidal closed category. For various questions in ``quantum'' algebra (see \cite{D}), we want to
know when the functor \[ \exists_N : [\A,\V] \too [\C,\V] \] \[ (\mbox{or } \Lan_N : \power (\A) \too \power (\C)) \]
given by the coend formula \[ \exists_N(f) = \int^a f(a) \tens \C(Na,-) \] (see \cite{DK} for notation),
is conservative (that is, reflects isomorphisms) for a given $\V$-functor \[ N : \A \too \C \] with $\A$ a small $\V$-category. In
this note we establish a ``simple'' sufficient condition for this to hold; namely, that $\A$ should also have a natural $\V$-opcategory structure
(with mild assumptions on $\V$).

\section{The main result}

If a small $\V$-category $\A$ is equipped with $\V$-natural transformations \begin{align*}
	\delta &: \A(a,b) \too \A(c,b) \tens \A(a,c) \\
	\epsilon &: \A(a,a) \too I
\end{align*} such that \[ \bfig \Square/>`>`>`<-/[\A(a,b)`\A(b,b) \otimes \A(a,b)`\A(a,b)`I \tens \A(a,b);\delta`1`\epsilon \tens 1`\cong] \efig \]
commutes, 
then the $\V$-functor \[ \exists_N = \Lan_N : [\A,\V] \too [\C,\V] \] is conservative for a given functor $N : \A \too \V$ if we suppose that each
component \[ N : \A(a,b) \too \C(Na,Nb) \] is a regular mono (that is, the kernel of some pair of maps) and that the composite of regular monos
in $\V$ is again a regular mono, and the functor \[ - \tens X : \V \too \V \] preserves regular monos for each $X$ in $\V$.

The rest of this section is concerned with the proof of this statement.

Note: We shall we use the term ``opcat'' to refer to any application using a $\V$-natural transformation of the form \[ fb \too \int_x \A(x,b) \tens fx \]
derived from the $\V$-natural transformation $\delta$ by use of the Yoneda expansion of a given
$\V$-functor $f$. In fact, we can suppose that $\A$ is merely a Frobenius category in the sense that the given family of maps
\[ \delta: \A(a,b) \too \A(c,b) \tens \A(a,c) \] is only $\V$-natural in $a$ and $b$ (and not necessarily in $c$) and that the family
\[ \epsilon : \A(a,a) \too I \]
is not necessarily $\V$-natural in $a$. We then use the $\V$-natural transformation \[ fb \too \prod_x \A(x,b) \tens fx \] in the following
calculations, where $\prod_x$ replaces $\int_x$, etc.

To show that $\exists_N$ is conservative, it suffices to show that the unit of the $\V$-adjunction $\exists_N \ladj [N,1]$ is a regular mono (see
\cite{BW} or \cite{KP}, for example, and the references therein to W. Tholen). To do this, we now consider the following diagram in $\V$.
\[ \bfig

%\!\!\!\!\!\!\!\!\!\!\!\!
%\!\!\!\!\!\!\!\!\!\!\!\!
%\!\!\!\!\!\!\!\!\!\!\!\!
\!\!\!\!\!\!\!\!\!\!\!\!
\!\!\!\!\!\!\!\!\!\!\!\!
\hscalefactor{1.8}
\node a(0,+1000)[\int^a \A(a,b) \tens fa]
\node b(+1000,+1000)[\int^a \C(Na,Nb) \tens fa]
\node c(0,+500)[\int^a \int_x \A(a,x) \tens \A(x,b) \tens fa]
\node d(+1000,+500)[\int^a\int_x \A(a,x) \tens \C(Nx,Nb) \tens fa]
\node e(0,0)[\int_x\int^a \A(a,x) \tens \A(x,b) \tens fa]
\node f(+1000,0)[\int_x\int^a \A(a,x) \tens \C(Nx,Nb) \tens fa]
\node h(-500,-500)[\int_x \A(x,b) \tens fx]
\node i(+1000,-500)[\int_x \C(Nx,Nb) \tens fx]
\node g(-500,+500)[fb]

\arrow|a|[a`b;\int^a N \tens 1]
\arrow|m|[a`c;\rm opcat]
\arrow|m|[b`d;\rm opcat]
\arrow[c`d;]
\arrow|m|[c`e;\rm can]
\arrow|m|[d`f;\rm can]
\arrow[e`f;]
\arrow|m|[e`h;\cong]
\arrow|m|[f`i;\cong]
\arrow|b|[h`i;\int_x N \tens 1]
\arrow|a|[a`g;\cong]
\arrow|m|[g`h;\rm opcat]

\efig \] where the isomorphisms are by the Yoneda expansion of $f$, and ``can'' denotes the canonical ``interchange'' maps. Since $\int_x N \tens 1$
is a regular mono (by the hypotheses on $\V$ and $N$), so also is its pullback along the right-hand composite map. Hence, since \[ {\rm opcat} : fb \too
\int_x \A(x,b) \tens fx \] is a coretraction, with left inverse the composite 
\[ \bfig \Square/<--`<-`>`<-/[fb`\int_x \A(x,b) \tens fx`I \tens fb`\A(b,b) \tens fb;%
`\cong`{\rm proj}_{x=b}`\epsilon \tens 1] \efig, \] we have that \[ \int^a N \tens 1 \] is the composite of a coretraction (which is always a regular mono)
and a regular mono (the pullback of $\int_x N \tens 1$), so is a regular mono (by hypothesis on $\V$). Thus, the adjunction unit we are looking at,
which is the composite \[ fb \oto{\cong} \int^a \A(a,b) \tens fa \oto{\int^a N \tens 1} \int^a \C(Na,Nb) \tens fa \] is a regular monomorphism, as
required.

\section{A frequent generalisation}

Similarly, given a $\V$-functor $N : \A \too \E$ with $\A$ small and $\E$ a $\V$-cocomplete $\V$-category, we have the 
standard $\V$-adjunction \[ (\epsilon,\eta):
L \ladj Y : \E \too \left[ \A^{\rm op},\V \right] \] where $Y(e)(a) = \E(Na,e)$ describes the $N$-Yoneda functor $Y$, and \[ L(f) = \int^a fa \tens Na\]
describes its left adjoint.
Hence, with corresponding hypotheses on $\A$, $N$, and $\V$ as before, we can replace the right-hand side of the diagram displayed in \S2 by
the following composite (*): \[ \bfig

\hscalefactor{1.8}
\node a(-500,+500)[\displaystyle \E\left(Nb,\int^a fa \tens Na \right)]
\node b(+500,+500)[\displaystyle\E\left( Nb, \int^a \left( \int_x \A(a,x) \tens fx \right) \tens Na \right)]
\node c(+500,0)[\displaystyle\E\left( Nb, \int^a \int_x \left( \A(a,x) \tens fx \right) \tens Na \right)]
\node e(+500,-500)[\displaystyle\E\left( Nb, \int_x \int^a \left( \A(a,x) \tens fx \right) \tens Na \right)]
\node d(-500,-500)[\displaystyle\E\left(Nb,\int_x fx \tens Nx \right)]

\arrow|l|/-->/[a`d;\rm (*)]
\arrow|a|[a`b;\rm opcat]
\arrow|r|[b`c;\rm can]
\arrow|r|[c`e;\rm can]
\arrow|b|[e`d;\cong]

\efig \]
\noindent where the isomorphism comes from the Yoneda-lemma expansion 
\[ Nx \to/<-/^{\cong} \int^a \A(a,x) \tens Na \] for $N$. Then, by the same argument as before,
we obtain a regular mono in $\V$ from the commuting diagram \[ \bfig \Square[fb`\displaystyle \E\left(Nb,\int^a fa \tens Na\right) = YL(f)(b)`%
\displaystyle \int_a fa \tens \A(b,a)`\displaystyle \E\left( Nb, \int_a fa \tens Na\right);\eta_{f,b}`{\rm opcat}`{\rm (*)}`] \efig \]
\noindent for each unit-component $\eta_{f,b}$ of the $\V$-functor $f$ in $[\A^{\rm op},\V]$, provided the canonical natural transformation
\[ X \tens \E(Nb,Na) \too \E(Nb, X \tens Na) \tag{$\dagger$} \]
is a regular mono in $\V$ for each $X$ in $\V$. Consequently, if this additional condition $(\dagger)$ holds
on $N$ and $\V$, then the left-adjoint $\V$-functor \[ L : [\A^{\rm op},\V] \too \E \] is conservative.

The result of \S2 can then be recovered from this latter result by putting $\E = [\C,\V]$, and taking the new $N: \A \too \E$ to be the composite of the 
$\V$-functor $N^{\rm \op} : \A^{\rm op} \too \C^{\rm op}$ (this $N$ from \S2) with the Yoneda embedding \[ \C^{\rm op} \subset [\C,\V], \]
noting that $\A$ is a $\V$-opcategory iff $\A^{\rm op}$ is also a $\V$-opcategory. The condition $(\dagger)$ holds for the new $N$ since
we have the isomorphisms

\[ \bfig \Square/>`<-`<-`<->/[X \tens \lbrak\C,\V\rbrak\left( \C(Nb,-),\C(Na,-)\right)`\lbrak\C,\V\rbrak\left( \C(Nb,-),X \tens \C(Na,-)\right)`%
X \tens \C(Na,Nb)`X \tens \C(Na,Nb);\cong`\cong`\cong`=] \efig \] by the Yoneda lemma applied twice.

\section{Related Conditions}

The following result is related to that of the earlier \S 2, but is much simpler.

Suppose that regular monomorphisms (that is, kernels in $\V$) are closed under composition in $\V$ and also that $n$ is a regular mono in $\V$ if
$mn$ and $m$ are regular monos in $\V$. Then, provided both coproduct $\Sigma$ and tensoring $X \tens -$ preserve regular monos in $\V$, we have that
\[ \exists_N : [\A,\V] \too [\C,\V] \]
\noindent is conservative if regular epimorphisms (that is, cokernels) split in the functor category $[\A^{\rm op} \tens \A,\V]$ and each
component \[ N : \A(a,b) \too \C(Na,Nb) \] of $N : \A \too \C$ is a regular monomorphism in $\V$.

To establish this, one simply notes that the canonical regular epimorphism \[ \sum_a S(a,a) \too \int^a S(a,a) \] splits naturally
in $S \in [\A^{\rm op} \tens \A,\V]$, so that $\int^a N \tens 1$ is a regular monomorphism by the hypotheses on $N$ and $\V$. This implies
that the unit components \[ \eta_{f,b} = fb \oto{\cong} \int^a \A(a,b) \tens fa \oto{\int^a N \tens 1} \int^a \C(Na,Nb) \tens fa \]
are regular monomorphisms, as required for $\exists_N$ to be conservative.

Note: In practice, it is often the case that, under the given hypotheses on $[\A^\op \tens \A,\V]$, each monomorphism of the form
\[ N : \A(-,b) \too \C(N-,Nb) \] splits naturally in $[\A^\op,\V]$, in which case the result is obvious and generalises the familiar fact that
$\exists_N$ is fully faithful if $N$ is.

Neither the condition in \S2 and \S3 that $\A$ should have an opcategory structure, nor the condition in this section that regular epimorphisms
should split in $[\A^\op \tens \A,\V]$, is in any way necessary for the process of left Kan extension along the Yoneda embedding of $\A^\op$
into $[\A,\V]$, to be conservative. One notable example is the (left) Cayley functor \[ \exists_P : [\A,\V] \too [\A^\op \tens \A,\V], \]
which is given by the coend formula \[ \exists_P(f) = \int^a f(a) \tens P(a,-,-), \] and is not fully faithful in general. This functor is
conservative for each (small) $\V$-promonoidal category $(\A,P,J)$ defined over \emph{any} complete and cocomplete base category $\V$.

Finally we note that both the sufficient conditions mentioned immediately above are closely related to the splitting properties of regular
epimorphisms in the functor category $[\A,\V]$. (by Maschke's result).

Note: Any enquires regarding this article can be forwarded to the author through Micah McCurdy (Macquarie University), who kindly typed the
manuscript.

\end{document}